\documentclass[12pt]{article}\usepackage{amsfonts,amssymb}

\usepackage{amsmath,bbm}
\usepackage{theorem}          
\usepackage{latexsym}
\usepackage{float}                   
\usepackage{times}
\usepackage{cite}

\newtheorem{tht}{Theorem}[section]
\newtheorem{thl}[tht]{Lemma}

\theoremstyle{plain}
{\theorembodyfont{\rmfamily}

\newtheorem{thex}{Example}
}

\newcommand{\mn}{\medskip}    

\newcommand{\sn}{\smallskip\noindent}

\newcommand{\cE}{{\mathcal{E}}} 
\newcommand{\cS}{{\mathcal{S}}}
\newcommand{\cC}{{\mathcal{C}}}
\newcommand{\cX}{{\mathcal{X}}}
\newcommand{\cD}{{\mathcal{D}}} 
\newcommand{\Hh}{{\mathcal{H}}} 
\newcommand{\cK}{{\mathcal{K}}}
\newcommand{\cU}{{\mathcal{U}}}  
\newcommand{\cT}{{\mathcal{T}}}
\newcommand{\cY}{\mathcal{Y}}
\newcommand{\cA}{\mathcal{A}}

\newcommand{\dR}{\mathbb{R}}
\newcommand{\dN}{\mathbb{N}}
\newcommand{\dC}{\mathbb{C}}
\newcommand{\dZ}{\mathbb{Z}}

\newcommand{\g}{\mathfrak{g}}
\newcommand{\gs}{{\mathfrak{s}}}

\newcommand{\im}{\mathrm{i}}
\newcommand{\I}{{\mathrm{Im}}}

\newcommand{\ovU}{\overline{dU(a)}}
\newcommand{\ovV}{\overline{dV(x_k)}}
\newcommand{\ov}{\overline}

\begin{document}

\date{\small{Fakult\"at f\"ur Mathematik und 
Informatik\\ Universit\"at Leipzig, 
Augustusplatz 10, 04109 Leipzig, Germany\\ 
E-mail: schmuedgen@math.uni-leipzig.de }
}

\title{A Strict Positivstellensatz for Enveloping Algebras}

\author{Konrad Schm\"udgen}

\maketitle

\renewcommand{\theenumi}{\roman{enumi}}
\begin{abstract}
Let $G$ be a connected and simply connected real Lie group with Lie algebra $\g$. Semialgebraic subsets of the unitary dual of $G$ are defined and a strict Positivstellensatz for positive elements of the universal enveloping algebra $\cE (\g)$ of $\g$ is proved.
\end{abstract}
{\small AMS subject classification: Primary: 14 P 10, 17 B 35}
\section{Introduction}
Positive polynomials on semialgebraic sets have been  intensively studied since E. Artin's solution of Hilbert's 17-th problem. In the last decade a number of new representation theorems for positive polynomials, usually called "Positivstellens\"atze", have been found (see e.g. \cite{S1}, \cite{P}, \cite{R}, \cite{PV}). Excellent surveys are given in the recent books \cite{PD}, \cite{M} and  the article \cite{S}. In a previous paper \cite{S3} a variant of  a non-commutative Positivstellensatz for the Weyl algebra was obtained. The aim of this paper is to prove a strict Positivstellensatz for enveloping algebras of finite dimensional Lie algebras.

Let $G$ be a connected and simply connected real Lie group with Lie algebra $\g$ and let $\cE(\g)$ be the complex universal enveloping algebra of $\g$. The algebra $\cE(\g)$ is a $\ast$-algebra with involution determined by $x^\ast = -x$ for $ x \in \g$. Let $\mathsf{1}$ denote the unit element of  $\cE(\g)$. Let $\{ x_1, {\dots}, x_d\}$ be a basis of $\g$ which will be fixed throughout this paper.

The algebra $\cE(\g)$ has a canonical filtration $(\cE_n(\g))_{n \ge 0}$, where $\cE_0(\g) = \dC \cdot \mathsf{1}$ and $\cE_n (\g), n \ge 1$, is the linear span of $\mathsf{1}$ and products $z_1, {\dots}, z_r$ with $z_1, {\dots}, z_r \in \g$ and $r \le n$ (see e.g. \cite{D}, 2.3). The associated graded algebra is the polynomial algebra $\dC [t_1,{\dots}, t_d]$, where the monomial $t^{k_1}_1 {\cdots} t^{k_d}_d$ corresponds to the element $x^{k_1}_1 {\cdots} x^{kd}_d$ of $\cE (\g)$. For an element $c \in \cE (\g)$  of degree $n$, we denote by $c_n (t)$ the polynomial of $\dC [t_1, {\dots}, t_d]$ corresponding to the component of $c$ with degree $n$.

Set $x_0 := \im\cdot 1$, where $\im$ denotes the complex unit. Then we have $x^\ast_j = - x_j$ for $j = 0, {\dots}, d$. Define
\begin{equation}\label{adef}
a := x^\ast_0 x_0 + x^\ast_1 x_1+ {\cdots}+ x^\ast_d x_d = \mathsf{1} - x^2_1 - {\cdots} - x^2_d.
\end{equation}
Let $S$ be a right Ore subset of $\cE(\g)\backslash \{0\}$ containing $a$. That is,  for any $s \in S$ and $z \in \cE (\g)$ there are elements $s^\prime \in S$ and $z^\prime \in \cE (\g)$ such that $s z^\prime = z s^\prime$. Note that such that a set $S$ exists since $\cE (\g)\backslash \{0\}$ is a (left and right) Ore set of $\cE (\g)$ (see \cite{D}, 3.6). For instance, if $a$ belongs to the center of $\cE(\g)$, then we may take the set of elements $a^n$, $n \in \dN_0$, as $S$.

Let $f = (f_1, {\dots}, f_r)$ be a finite set of hermitean elements of the enveloping algebra $\cE (\g)$ such that $f_1 = \mathsf{1}$. Let $\cK_f$ and  $\cT_f$ be the associated basic closed semialgebraic set and positive wedge, respectively, as defined by formulas (\ref{kset}) and (\ref{tfdef}) below.

The main result of this paper is the following
\begin{tht}\label{T1}
Suppose that $c$ is a hermitean element of  the enveloping algebra $\cE (\g)$ of even degree $2m$ satisfying the following assumptions:\\
(i) There exists $\varepsilon > 0$ such that $c - \varepsilon \cdot \mathsf{1} \in \cT_f$.\\
(ii) $c_{2m} (t) > 0$ for all $t \in \dR^d, t \ne 0$.\\
If $m$ is even, there exists an element $s \in S$ such that $s^\ast c s \in \cT_f$. If $m$ is odd, there is an $s \in S$ such that
$$
\sum^d_{k=0} s^\ast x^\ast_k c x_k s \in \cT_f.
$$
\end{tht}
This paper is organized as follows.
In Section 2 we collect some basics on the infinitesimal representation $dU$ of a unitary representation $U$ of $G$ and we define semialgebraic subsets of the unitary dual of $G$. In Section 3 we introduce an auxiliary $\ast$-algebra of bounded operators associated with the (unbounded) $\ast$-algebra $dU(\cE (\g))$. Representations of this algebra are studied in Section \ref{S4}. These results are used in the proof of Theorem 1.1 given in Section \ref{S5}.

Let $\cA$ be a unital complex $\ast$-algebra. We denote by $\sum^2 (\cA)$ the set of all finite sums of squares $x^\ast x$, where $x \in \cA$. A subset of the hermitean part $\cA_h := \{ x \in \cA : x^\ast =x\}$ is called an $m$-admissible wedge (\cite{S2}, p.22) if $\cC$ is a wedge (that is, $x + y \in \cC$ and $\lambda x \in \cC$ if $x, y \in \cC$ and  $\lambda \ge 0 \}$ such that the unit element is in $\cC$ and $z^\ast x z \in \cC$ for $x \in \cC$ and $z \in \cA$.
\section{Unitary Representations and Semialgebraic Sets of the Dual }
By a unitary representation of the Lie group $G$ we mean a strongly continuous homomorphism $U$ of $G$ into the group of unitary operators of a Hilbert space $\Hh (U)$. Let $\cD^\infty (U)$ denote the vector space of $C^\infty$-vectors of $U$ and let $dU$ be the associated $\ast$-representation of the $\ast$-algebra $\cE (\g)$ on the dense domain $\cD^\infty (U)$ of $\Hh(U)$ (see \cite{S2}, Chapter 10, or \cite{Wa}, Section 4.4, for details). For $f \in \cE (\g)$ we write $dU(f) \ge 0$ when $\langle d U(f) \varphi, \varphi \rangle \ge 0$ for all vectors $\varphi \in \cD^\infty (U)$.

For later use we restate some classical results of E. Nelson and W. F. Stinespring \cite{NS} and  of E. Nelson \cite{N} in 
\begin{thl}\label{L2}
If $U$ is a unitary representation of $G$, then the closure $\ov{d U (a)}$ of the operator $d U (a)$ is self-adjoint and equal to  $B := I - \sum^d_{k=1} \ov{d U (x_k)}^2$. Moreover,
\begin{equation}\label{dinfvec}
\cD^\infty (U) = \cap^\infty_{n=1} \cD (\ovU^n) = \cap^\infty_{n=1} \cD (B^n).
\end{equation}
\end{thl}
{\bf Proof.} For the self-adjointness of $\ovU$ and  the first equality of (\ref{dinfvec}) see Corollaries 10.2.4 and 10.2.7 in \cite{S2} or Theorems 4.4.3 and 4.4.4 in \cite{Wa}. We prove that $\ovU = B$. By Lemma 10.4.5 in \cite{S2} or by Lemma 4.4.4.8 in \cite{Wa} there is constant $c > 0$ such that
$$
\| d U (x_k) \varphi \| + \| dU(x_k)^2 \varphi \| \le c \| d U (a) \varphi \|, \varphi \in\cD^\infty (U), k = 1, {\dots}, d.$$

Since $B \varphi = d U (a) \varphi$ for $\varphi \in \cD^\infty (U)$, the preceding inequality implies that $\cD (\ovU) \subseteq \cD (B)$ and $B \varphi = \ovU \varphi$ for $\varphi \in \cD (\ovU)$.
Since $B$ is a symmetric extension of the self-adjoint operator  $\ovU$, we get $B = \ovU$. \hfill $\Box$

\sn
Let $\hat{G}$  denote the unitary dual of $G$, that is, $\hat{G}$ is the set of unitary equivalence classes of irreducible unitary  representations of $G$. For each $\alpha \in \hat{G}$ we fix a representation $U_\alpha$ of the equivalence class $\alpha$.

\sn
{\bf Definition:} A subset $K$ of $\hat{G}$ is called {\it semialgebraic} if $K$ is a finite Boolean combination (that is, using unions, intersections and complements) of sets $K_f = \{ \alpha \in \hat{G} : d U_\alpha (f) \ge 0\}$ with $f \in \cE (\g)$.

\sn
For $r\in \dN$ and an $r$-tuple $f = (f_1, {\dots}, f_r)$ of elements $f_j \in \cE (\g)$ such that $f_1 = 1$, we define the basic closed semialgebraic set
\begin{equation}\label{kset}
\cK_f = \{ \alpha \in \hat{G} : d U_\alpha (f_1) \ge 0,{\dots}, d U_\alpha (f_r) \ge 0\}
\end{equation}
and  the associated positive wedges
\begin{equation}\label{tfdef}
\cT_f = \{ z = \sum^k_{j=1} \sum^r_{l=1} z^\ast_{jl} f_l z_{jl}; z_{jl} \in \cE (\g), k \in \dN\},
\end{equation}
$$
\cE(\g; f)_+ = \{z \in \cE(\g)_h : dU_\alpha (z) \ge 0~~{\rm for}~~\alpha \in K_f).
$$
Clearly, $\cT$ and $\cE (\g;f)_+$ are $m$-admissible wedges of $\cE(\g)$ such that $\cT_f \supseteq \cE(g;f)_+$.

If $f = (1)$, then $\cK_f = \hat{G}$. In this case we denote $\cE (\g;f)_+$ by $\cE(\g)_+$. From decomposition theory (see e.g. \cite{S2}) it follows that $\cE(\g)_+$ is the set of elements $x \in \cE(\g)$ such that $dU(x)\ge 0$ for all unitary representations $U$ of the group $G$.
\begin{thex}\label{Ex1}
$G=\dR^d$\\
Then $g$ is the abelian Lie algebra $\dR^d$ and $\cE(\g)$ is the  polynomial algebra \break $\dC[x_1,{\dots}, x_d]$. For each point $t=(t_1,{\dots},t_d) \in\dR^d$ there is an irreducible unitary representation $U_t$ of the Lie group $\dR^d$ on $\Hh(U_t)=\dC$ such that $dU_t(x_j)=t_j,j=1,{\dots},d$. Because  these representations $U_t$ exhaust the dual of $\dR^d$, we can identify the dual of the Lie group $\dR^d$ with $\dR^d$. In this manner the semialgebraic sets according to the definition given above are just the "ordinary" semialgebraic sets in semialgebraic geometry (\cite{PD}, p. 31). For an $r$-tuple $f=(f_1,{\dots}, f_r)$  elements $f_j\in \dC[x_1,{\dots}, x_d]$ with $f_1=1$, let $\tilde{f}$ be the $\binom{r}{2} +1$-tuple of elements $f_1$ and $f_i f_j,i\ne j, i,j=1,{\dots},r$. Then we have $\cK_f=\cK_{\tilde{f}}$. The wedge $\cT_{\tilde{f}}$ is the usual preorder from semialgebraic geometry  and $\cE(\g;f)_+$ is the set of real polynomials which are nonnegative on $\cK_f$. Note that for noncommutative Lie groups products $f_if_j$ cannot be added to the wedge $\cT_f$ in general because the product of two noncommuting hermitean elements of $\cE(\g)$ is not even hermitean.
\end{thex}
\begin{thex}\label{Ex2} $G=SU(2)$\\
The Lie algebra of $SU(2)$ has a basis $\{x_1,x_2,x_3\}$ satisfying the commutation relations 
\begin{equation}\label{1xrel}
[x_1, x_2]=x_3, [x_2,x_3]= x_1,[x_3,x_1]=x_2.
\end{equation} 
The element $a = 1 - x^2_1 - x^2_2 - x^2_3$ generates the center of the algebra $\cE(\g)$. The unitary dual of $G$ consists of equivalence classes of spin $l$ representations $U_l$, $l \in \frac{1}{2}\dN_0$. We identify $\hat{G}$ and $\frac{1}{2} \dN_0$. Set $H = -i x_1$. Here we only need that the operators $dU_l(H)$ and $dU_l(a)$ act on  an orthonormal basis $e^l_j, j=-l,-l+1,{\dots},l$, of a $(2l+1)$-dimensional Hilbert space by
\begin{equation*}
dU_l (H) e_j^l =2j e^l_{j},~ dU_l(a) e^l_j = (l^2+l+1) e^l_j.
\end{equation*}

Clearly, a polynomial $p(H)$ in the generator $H$ is in $\cE(\g)_+$ if and only if $p(n)\ge 0$ for all $n\in\dZ$. Using the relations (\ref{1xrel}) it is easily shown that $p(H)\in \sum^2 (\cE(\g))$ if and only if there is a polynomial $q(H)\in \dC[H]$ such that $p(H)=q(H)^\ast q(H)$. In particular, $(H-c_1)(H-c_2) \in \cE(\g)_+$ and $(H-c_1)(H-c_2)\not\in \sum^2 (\cE(\g))$ if there is in integer $n$ such that $n\le c_1 <c_2\le n+1$.

Let $M=\{l_1,{\dots},l_r\}$ be a finite subset of $\frac{1}{2} \dN_0 \cong \hat{G}$. Put $f_1 = g_1=\mathsf{1} , f_{j+1} = 2l_j-H$ and $g_{j+1}=-(a-l_j^2-l_j-1))^2$ for $j=1, \dots,r$. Then $\cK_f = \cK_g=M$, so $M$ is a semi-algebraic subset o f $\hat{G}$.
\end{thex}

\section{An Auxiliary $\ast$-Algebra}\label{S3}
In what follows $U$ denotes a fixed unitary representation of the Lie group $G$ on a Hilbert space $\Hh (U)$. In this section we define and study an auxiliary $\ast$-algebra $\cX$  depending on $U$.

\mn 
For notational simplicity we abbreviate 
$$
X_k := dU (x_k), k = 0, {\dots}, d,~~{\rm and}~~ A := dU(a).
$$ 
Clearly, $\langle A \varphi, \varphi\rangle \ge \langle \varphi,\varphi\rangle$ for $\varphi \in \cD^\infty (U)$. Hence the inverse $Y := A^{-1}$ exists and  maps $\cD^\infty (U)$ onto $\cD^\infty (U)$ by Lemma \ref{L2}.

Let $\cX$ denote the algebra of operators acting on the invariant dense domain $\cD^\infty (U)$ of $\Hh (U)$ generated by the identity map $I \equiv I_\cD$, 
\begin{equation}\label{ydef}
Y_{kl} := X_k X_l Y~~{\rm for}~~ Y_{-l,-k} := Y X_l X_k~~{\rm and}~~ k, l = 0,{\dots}, d.
\end{equation}
Clearly, $\cX$ is a $\ast$-algebra of operators with involution
\begin{equation}\label{yadj}
Y^\ast_{kl} = Y_{-l,-k}, k, l = 0,{\dots},d.
\end{equation}
Let $\cX_0$ denote two-sided $\ast$-ideal of $\cX$ generated by $Y_{k0}, k = -d, {\dots}, d$. Let $c^k_{lj}$ denote the structure constants of the Lie algebra $\g$ and set $b^k_{ij} := c^k_{ij} + c^j_{ik}$. Then we have
\begin{equation}\label{lie}
X_i X_j - X_j X_i = \sum^d_{k=1} c^k_{ij} X_k.
\end{equation}
The operators $Y_{kl}, X_j$ and $Y = -Y_{00}$ satisfy the following relations:
\begin{align}\label{r1}
&\sum^d_{k=0} Y^\ast_{k0} Y_{k0} = Y,\\
\label{r2}
&Y^\ast_{k0} Y_{l0} - Y^\ast_{l0} Y_{k0} = -\im \sum^d_{j=1}  c^j_{lk} Y Y_{j0}~~{\rm for} ~~ k, l = 1,{\dots},d,\\
\label{r3}
&Y^\ast_{k0} - Y_{k0} = \im \sum^d_{j,l=1} b^j_{kl} Y^\ast_{j0} Y_{l0}~~{\rm for}~~ k = 1, {\dots}, d,\\
\label{r4}
&\sum^d_{k,l=0} Y^\ast_{kl} Y_{kl} = I + \im \sum^d_{j,k,l=1} b^l_{jk} Y^\ast_{j0} Y_{kl} = I - \im \sum^d_{j,k,l =1} b^l_{jk} Y^\ast_{kl} Y_{j0},\\
\label{r5}
&Y_{kl} - Y^\ast_{kl} \in \cX_0, Y_{kl} - Y_{lk} \in \cX_0~~{\rm for}~~ k,l = 0,{\dots},d,\\
\label{r6}
&Y x - x Y \in Y \cX_0 = \cX_0 Y~~{\rm for}~~ x \in \cX,\\
\label{r7}
&Y_{kl} Y_{ij} - Y_{ij} Y_{kl} \in \cX_0~~{\rm for}~~  k, l, i, j = 1,{\dots}, d,\\
\label{r8}
&Y_{kl} Y_{ij} - Y_{ki} Y_{lj} \in \cX_0~~{\rm for}~~  k, l, i, j = 1,{\dots}, d.
\end{align}
Note that relations (\ref{r1}) -- (\ref{r3}) above have been found by W. Szymanski \cite{Sz}. They are relations (3) -- (5) in \cite{Sz}. Setting $ l = 0$ in (\ref{yadj}) we get relations (1) -- (2) in \cite{Sz}.

The relations (\ref{r1}) -- (\ref{r8}) listed above are proved by straightforward algebraic manipulations using  the Lie algebra commutation relation (\ref{lie}). Using (\ref{adef}), (\ref{lie}) and the  abbreviations $b^k_{ij} = c^k_{ij} + c^j_{ik}, Y = d U(a)^{-1}$ we obtain
\begin{align}\label{yxk1}
&Y X_k = X_k Y + \sum^d_{i,j=1} b^j_{kl} Y X_l X_j Y,\\
\label{yxk2}
&Y X_k X_l = X_k X_l Y + \sum^d_{i,j=1} (b^j_{li} Y X_k X_i X_j Y + b^j_{ki} Y X_i X_j X_l Y)
\end{align}
for $k,l = 1, {\dots}, d$. All relations (\ref{r1}) -- (\ref{r8}) are easily derived from (\ref{yxk1}) and (\ref{yxk2}) combined with (\ref{ydef}). We omit the details of these verifications.
\begin{thl}\label{L3}
For arbitrary numbers $n \in \dN$ and $k_1, {\dots}, k_{4n} \in \{0, {\dots}, d\}$, we have 
\begin{align}\label{yx1}
&Y^n X_{k_1} {\cdots} X_{k_{2n}} \in \cX\\
\label{yx2}
&Y^n X_{k_1} {\cdots} Y_{k_{4n}} Y^n - Y_{k_1 k_2} {\cdots} Y_{k_{4n-1}, k_{4n}} \in \cX_0.
\end{align}
\end{thl}
{\bf Proof.} Both equations are proved by induction on $n$. We carry out the proof of (\ref{yx2}). First let $n =1$. By (\ref{ydef}) and (\ref{r5}) we get 
$$
Y X_{k_1} X_{k_2} X_{k_3} X_{k_4} Y - Y_{k_1 k_2} Y_{k_3 k_4} = (Y^\ast_{k_2 k_1}  - Y_{k_1 k_2}) Y_{k_3 k_4}\in \cX_0.
$$
Suppose the assertion holds for $n$. Using the abbreviations 
$z_n = Y^n X_{k_1} {\cdots} X_{k_{4n}} Y^n$, $z = Y X_{k_{4n+1}} {\cdots} X_{k_{4n+4}} Y, y_n = Y_{k_1 k_2} {\cdots} Y_{k_{4n-1}, k_{4n}}, y = Y_{k_{4n+1}, k_{4n+2}} Y_{k_{4n+3}, k_{4n+4}}, Y = A^{-1}$ we compute
\begin{align}\label{right}
&Y^{n+1} X_{k_1} {\cdots} X_{k_{4n+4}} Y^{n+1} -  Y_{k_1 k_2} {\cdots} Y_{k_{4n+3},k_{4n+4}} = Y z_n A^{n+1} z Y^n - y_n y \hspace{3cm}\nonumber\\
&= [Y,z_n] A^{n+1} [z, Y^n] + z_n A^n [z, Y^n] + [Y, z_n] Az + (z_n {-} y_n) z + y_n (z {-} y).
\end{align}
From (\ref{r6}) we easily derive that $[z, Y^n] \in Y^n \cX_0$.
Moreover, $[Y, z_n] \in Y \cX_0 = \cX_0 Y$ by (\ref{r6}) and $z_n - y_n \in \cX_0$ and $z - y \in \cX_0$ by the induction hypothesis. Using these facts and remembering that $\cX_0$ is a two-sided ideal of $\cX$ it follows that the element in (\ref{right}) belongs to $\cX_0$. This proves the assertion for $n+1$.\hfill $\Box$

\sn
For $z \in \cX$ we set $Re\ z := \frac{1}{2} (z + z^\ast)$ and $\I \  z := \frac{1}{2} \im (z^\ast - z)$. Let $\cX_b$ be the set of all elements $z \in \cX$ for which there exists a positive number $\lambda$ such that
$$
\lambda \cdot I \pm ~{\rm Re}~ z \in \sum\nolimits^2 (\cX)~~{\rm and}~~ \lambda \cdot I \pm \I z \in \sum\nolimits^2 (\cX).
$$
From Lemma \ref{L2}(ii) in \cite{S3} it follows that a finite sum  $\sum_j z^\ast_j z_j$ is in $\cX_b$ if and only if all $z_j$ are in $\cX_b$. Moreover, $\cX_b$ is a $\ast$-algebra by Corollary 2.2 in \cite{S3}. We shall use these two facts in the proof of Lemma \ref{L4} below. Following \cite{S3} we say that the $\ast$-algebra $\cX$ is {\it algebraically bounded} if $\cX = \cX_b$.
\begin{thl} \label{L4}
The $\ast$-algebra $\cX$ defined above is algebraically bounded.
\end{thl}
{\bf Proof.}
Recall that $Y^\ast = Y = -Y_{00}$. Applying relation (\ref{r1}) we obtain 
$$
\mbox{$\frac{1}{4}$}  I - (\mbox{$ \frac{1}{2}$}  I - Y )^2 = Y - Y^2_{00}   =_{_{_{\!\!\!\!\! \!\!\! \!(r1)}}}  \sum^d_{k=1} Y^\ast_{k0} Y_{k0} \in \sum\nolimits^2 (\cX).
$$
Thus, $( \frac{1}{2} I - Y)^2 \in \cX_b$ and hence $Y \in \cX_b$. Since $Y \in \cX_b$, it follows from (\ref{r1}) that $Y_{k0} \in \cX_b$ for $k = 1,{\dots}, d$. Using relation (\ref{r4}) and the fact that $Y_{k0} = Y_{0k}$ we compute
\begin{align*}
&Y^2_{00} + 2 \sum^d_{k=1} Y^\ast_{k0} Y_{k0} + \sum^d_{k,l=1} (Y_{kl} + \frac{\im}{2} \sum^d_{j=1} b^k_{jl} Y_{j0})^\ast (Y_{kl} + \frac{\im}{2} \sum^d_{j=1} b^k_{jl} Y_{j0})\\
&= \sum^d_{k,l=0} Y^\ast_{kl} Y_{kl} - \frac{\im}{2} \sum^d_{j,k,l=1} b^k_{jl} Y^\ast_{j0} Y_{kl} + \frac{\im}{2} \sum^d_{j,k,l=1} b^k_{jl} Y^\ast_{kl} Y_{j0}\noindent\\
&~~~~+ \frac{1}{4}   \sum^d_{k,l=1} \left(  \sum^d_{j=1} b^k_{jl} Y_{j0} \right)^\ast \left(  \sum^d_{j=1} b^k_{jl} Y_{j0} \right)\\
& =_{_{_{\!\!\!\!\! \!\!\! \!(r4)}}} I + \frac{1}{4}   \sum^d_{k,l=1} \left(  \sum^d_{j=1} b^k_{jl} Y_{j0} \right)^\ast \left(  \sum^d_{j=1} b^k_{jl} Y_{j0} \right).
\end{align*}
Since $Y_{j0} \in \cX_b$ for $j = 0, {\dots}, d$ as just shown and $\cX_b$ is a $\ast$-algebra, the right-hand side of the preceding equation belongs to $\cX_b$. Therefore, $Y_{kl} + \frac{\im}{2} \sum^d_{j=1} b^k_{jl}  Y_{j0} \in \cX_b$ and hence $Y_{kl} \in \cX_b$ for $k,l = 1, {\dots}, d$. Hence all generators of $\cX$ are in $\cX_b$, so that $\cX = \cX_b$. \hfill $\Box$

\sn
Now we choose $p \in \dN$ such that $4p \ge$ degree $f_j$ for $j = 1, {\dots}, r$. Then, $Y^p dU (f_j) Y^p \in \cX$ by (\ref{yx2}). Let $\cC_f$ denote the wedge of all finite sums of elements $x^\ast x$ and $z^\ast Y^p dU (f_l) Y^p z$, where $x, z \in \cX$ and $ l=1,  {\dots}, r$. The assertion of the next lemma is contained in \cite{S3}, Lemma 2.3. For completeness we include the short proof.
\begin{thl}\label{S3-L5}
If $ z\in \cX$ is not in $\cC_f$, then there exists a state $F$ of the $\ast$-algebra $\cX$ such that $F(z)\le 0$ and $F(x)\ge 0$ for all $x\in \cC_f$.
\end{thl}
{\bf Proof.} Since $\cC_f \subseteq \sum^2 (\cX)$, the unit element $I$ of $\cX$ is an internal point of the wedge $\cC_f$ in the real vector space $\cX_h = \{ x \in \cX : x = x^\ast\}$. By the separation theorem for convex sets \cite{K}, \S 17, (3), there is a linear functional $G\not\equiv 0$ on $\cX_h$ such that $G (z) \le 0$ and $G(x)\ge 0$ for $x \in \cC_f$. Since $G \not\equiv 0$, we have $G (I) > 0$. Take as $F$ the extension of the $\dR$-linear functional $G(I)^{-1} G$ on $\cX_h$ to a $\dC$-linear functional on $\cX$. \hfill $\Box$
\section{Representations of the Auxiliary $\ast$-Algebra}\label{S4}
Since the $\ast$-algebra $\cX$ is  algebraically bounded by Lemma \ref{L4}, for any $\ast$-representa\-tion of $\cX$ all representation operators are bounded.
Let $\pi$ be an arbitrary $\ast$-representation of $\cX$ on a Hilbert space $\Hh$. By (\ref{r6}), $\Hh_\infty := \ker \pi (Y)$ is  invariant and hence reducing for the bounded $\ast$-representation $\pi$. Let $\pi_\infty$ and $\pi_0$ denote the restrictions of $\pi$ to $\Hh_\infty$ and $\Hh_0 := \Hh^\bot_\infty$, respectively.
\subsection {} \label{S4-sub1}
In this subsection we investigate the $\ast$-representation $\pi_0$. Since $\ker \pi_0 (Y) = \{0\}$ and relations (\ref{r1}) -- (\ref{r3}) hold, Lemma 1 in \cite{Sz} applies. It is a reformulation of Nelson's famous integrability theorem for Lie algebra representations (\cite{N}, see e.g. \cite{S2}, Theorem 10.5.6, or \cite{Wa}, Theorem 4.4.6.6) and states that there exists a unitary representation $V$ of the simply connected Lie group $G$ on $\Hh_0$ such that
\begin{equation}\label{dvx}
\ovV = -\im \ov{\pi_0 (Y_{k0}) \pi_0 (Y)^{-1}},~ k = 1, {\dots}, d.
\end{equation}
In the proof therein it is shown that $B:= I - \sum^d_{k=1} \ovV^2$ is equal to the self-adjoint operator $\pi_0 (Y)^{-1}$ on its domain $\pi_0 (Y) \Hh_0$. By Lemma 2.1, $\cD^\infty (V) = \cap^\infty_{n=1} \cD (B^n) = \cap^\infty_{n=1} \pi_0 (Y)^n \Hh_0$. Hence $\pi_0 (Y)^{-1}$ maps $\cD^\infty (V)$ onto $ \cD^\infty (V)$. Therefore, by (\ref{dvx}) we have
\begin{equation}\label{dvk}
 dV(x_k)\varphi = -\im \pi_0 (Y_{k0}) \pi_0 (Y)^{-1} \varphi,~ \varphi \in \cD^\infty (V), k = 1, {\dots}, d.
\end{equation}

\sn
Next we prove by induction on $n$ that 
\begin{equation}\label{pil}
\pi_0 (Y^n X_{k_1} {\cdots} X_{k_{2n}}) \varphi = \pi_0 (Y)^n dV(x_{k_1} {\cdots} x_{k_{2n}}) \varphi, ~\varphi \in \cD^\infty (V),
\end{equation}
for $k_1, {\dots}, k_{2n} \in \{0, {\dots}, d\}$. Since $Y^n  X_{k_1}{\cdots} X_{k_{2n}} \in \cX$ by (\ref{r7}), the left hand side of (\ref{pil}) is well-defined. First let $n = 1$. For $\varphi \in \cD^\infty (V)$ we set $\psi = \pi_0 (Y)^{-1} \varphi$. Using (\ref{ydef}) and (\ref{dvk}) we compute
\begin{align*}
&\pi_0 (Y X_{k_1} X_{k_2}) \varphi = \pi_0 (Y X_{k_1} X_{k_2}) \pi_0 (Y) \psi 
= \pi_0 (Y X_{k_1} Y) \pi_0 (Y)^{-1} \pi_0 (X_{k_2} Y) \psi \hspace{1cm} \nonumber \\
&= -\pi_0 (Y) \pi_0 (Y_{k_1, 0}) \pi_0 (Y)^{-1} \pi_0(Y_{k_2, 0}) \pi_0 (Y)^{-1} \varphi\nonumber\\
& = \pi_0 (Y) dV(x_{k_1}) dV(x_{k_2}) \varphi = \pi_0 (Y) dV (x_{k_1} x_{k_2}) \varphi
\end{align*}
which proves (\ref{pil}) for $n=1$. Suppose now that (\ref{pil}) is true for $n$. Let $\varphi \in \cD^\infty (V)$. Set $\psi = \pi_0 (V)^{-1} \varphi, z_n = X_{k_1} {\dots} X_{k_{2n}}$ and $z = X_{k_{2n+1}} X_{k_{2n+2}}$. Since $\pi_0$ and $dV$ are $\ast$-representations and $\pi_0 (Y z^\ast) = \pi_0 (Y) dV (z^\ast)$ on $\cD^\infty (V)$ (by (\ref{pil}) for $n = 1$), it follows that $\pi_0 (z Y) \psi = dV (z) \pi_0 (Y) \psi$. Using the latter relation and $\pi_0 (Y^n z_n) \psi = \pi_0 (Y)^n dV (z_n) \psi$ by the induction hypothesis we get
\begin{align*}
&\pi_0 (Y^{n+1} z_n z) \varphi = \pi_0 (Y) \pi_0 (Y^n z_n z Y) \psi 
= \pi_0 (Y) \pi_0 (Y^n z_n) \pi_0 (z Y) \psi \hspace{4cm}\\ 
&= \pi_0 (Y)^{n+1} dV (z_n) dV (z) \pi_0 (Y) \psi =\pi_0 (Y)^{n+1} dV(z_nz)\varphi
\end{align*}
which is equation (\ref{pil}) for $n + 1$. This completes the proof of (\ref{pil}). Applying the involution to both sides of (\ref{pil}) we obtain
\begin{equation}\label{pir}
\pi_0 (X_{k_1} {\cdots} X_{k_{2n}} Y^n) \varphi = dV (x_{k_1} {\cdots} x_{k_{2n}}) \pi_0 (Y)^n \varphi, \varphi \in \cD^\infty (V).
\end{equation}
Combining (\ref{pil}) and (\ref{pir}) we conclude that 
$$
\pi_0 (Y^n X_{k_1} {\cdots} X_{k_{4n}} Y^n) \varphi = \pi_0 (Y)^n dV (x_{k_1} {\cdots} x_{k_{4n}}) \pi_0 (Y)^n \varphi, \varphi \in \cD^\infty (V),
$$
for $k_1, {\dots} k_{4n} \in \{0, {\dots}, d\}$. This in turn implies that
\begin{equation}\label{pib}
\pi_0 (Y^n dU(x)Y^n) \varphi = \pi_0 (Y)^n dV (x) \pi_0 (Y)^n \varphi, \varphi \in \cD^\infty (V)\\
\end{equation}
for all $x \in \cE  (g)$ such that degree $x \le 4n$.
\subsection{}\label{S4-sub2}
In this subsection we turn to the $\ast$-representation $\pi_\infty$ of  $\cX$ on $\Hh_\infty = \ker \pi (Y)$. Since $\pi_\infty (Y) = 0$, it follows from (\ref{r1}) and  (\ref{yadj}) that
\begin{equation}\label{pinullinf}
\pi_\infty (Y_{k0}) = \pi_\infty (Y_{0k})~~{\rm for}~~k = -d, {\dots},d.
\end{equation}
Moreover, $\pi_\infty (\cX_0) = \{0\}$. Therefore, by (\ref{r5}), (\ref{r7}) and (\ref{r8}), $y_{kl} := \pi_\infty (Y_{kl})$, $k,l = 1, {\dots}, d$ and $ k,l = -d, {\dots}, -1$, are pairwise commuting bounded self-adjoint operators such that $y_{kl} = y_{lk}$ and $y_{-k,-l} = y_{kl}$. Let $\chi$ be a character of the abelian unital $C^\ast$-algebra $\cY$ generated by these operators (or equivalently by $\pi_\infty (\cX)$). From (\ref{r4}), (\ref{pinullinf}) and (\ref{r8}) we get 
\begin{equation}\label{ysum}
\sum^d_{k,l=1} y^2_{kl} =I~~{\rm and}~~ y_{kl} y_{ij} = y_{ki} y_{lj}, ~i,j,k,l = 1,{\dots}, d.
\end{equation}
From (\ref{ysum}) it follows that there is a $j \in \{1,{\dots}, d\}$ such that $\chi (y_{jj}) \ne 0$. Let $\epsilon$ denote the sign of $\chi (y_{jj})$. Take $t_j \in \dR$ such that $t^2_j = \epsilon \chi (y_{jj})$ and put \break
$t_k := \chi (y_{kj}) \chi (y_{jj})^{-1} t_j, k = 1, {\dots}, d$. By (\ref{ysum}) we have
\begin{equation}\label{chik}
\epsilon t_k t_l = \epsilon \chi (y_{kj}) \chi (y_{lj}) \chi(y_{jj})^{-2} t^2_j = \epsilon \chi (y_{kl} ) \chi (y_{jj})^{-1} t^2_j = \chi (y_{kl})
\end{equation}
for $k,l = 1,{\dots},d$ and hence
$$
\left( \sum^d_{k=1} t^2_k\right)^2 = \sum^d_{k,l=1} \epsilon \chi (y_{kk}) \epsilon \chi (y_{ll}) = \sum^d_{k,l =1} \chi (y^2_{kl}) = \chi (I) =1.
$$
Since all operators $y_{kl}$ are self-adjoint, all numbers $t_k$ are real. Thus, for each character $\chi$ of $\cY$ there exist $\epsilon \in \{-1,1\}$ and a point $t= ( t,{\dots}, t_d)$ of the unit sphere $S^d$ of $\dR^d$ such that (\ref{chik}) holds. From the Gelfand theory it follows that there are reducing subspaces $\Hh^\pm_\infty$ for $\pi_\infty$ such that $\Hh_\infty = \Hh^+_\infty \oplus \Hh^-_\infty$ and spectral measures $E^\pm$ over $S^d$ on $\Hh^\pm_\infty$ such that for $k,l=1,{\dots},d$,
\begin{equation}\label{piyinf}
\pi_\infty (y_{kl}) = \int_{S^d} t_k t_l dE^+ (t) \oplus - \int_{S^d} t_k t_l dE^-(t).
\end{equation}
\section{Proof of Theorem 1}\label{S5}
Suppose first that $m$ is even, say $m = 2n$. Since degree $c = 4n$, it follows from formula (\ref{yx2}) in Lemma \ref{L3} that $Y^n dU(c)Y^n \in \cX$.The crucial step of the proof is the assertion of the following
\begin{thl} \label{L6}
$Y^n dU(c)Y^n$ belongs to the wedge $\cC_f$ defined in Section \ref{S3}.
\end{thl}
{\bf  Proof.} Assume the contrary. Then,  by Lemma \ref{S3-L5}  there exists a state $F$ of $\cX$ such that $F(Y^n dU(c)Y^n)\le 0$ and $F(x) \ge 0$ for all $x \in \cC_f$. Let $\pi_F$ be the $\ast$-representation of $\cX$ associated with $F$ by the GNS construction. Then there is a cyclic vector $\varphi_F$ such that $F(x)=  \langle \pi_F (x) \varphi_F,\varphi_F\rangle, x \in \cX$. As  shown in Section \ref{S4}, $\pi_F$ decomposes into a direct sum of representations $\pi_0$ and $\pi_{\infty}$. If $\varphi_0$ and $\varphi_\infty$ are the corresponding components of $\varphi_F$, we have
\begin{equation}\label{fpix}
F(x) = \langle \pi_0 (x) \varphi_0,\varphi_0\rangle + \langle \pi_\infty (x) \varphi_\infty, \varphi_\infty\rangle,~ x\in \cX.
\end{equation}

Our next aim is to derive inequality (\ref{pic1}) below. Let $V$ be the unitary representation of $G$ from Subsection \ref{S4-sub1}.

We prove that $dV(f_l)\ge 0$ for $l = 1, {\dots}r$.  Let $\psi \in \cD^\infty (V)$. Since the $\ast$-representation $\pi_F$ is cyclic, there is a sequence $\{b_n; n \in \dN\}$ of elements $b_n \in \cX$ such that $\pi_0 (b_n) \varphi_0 \rightarrow \pi_0 (Y)^{-p} \psi$ and $\pi_\infty (b_n) \varphi_\infty \rightarrow 0$. From (\ref{fpix}) we obtain
\begin{align*}
& F(b^\ast_n Y^p d U(f_l) Y^p b_n)=\\
&\langle \pi_0 (Y^p dU(f_l) Y^p) \pi_0 (b_n) \varphi_0, \pi_0 (b_n) \varphi_0\rangle + \langle \pi_\infty (Y^p dU(f_l) Y^p)\pi_\infty (b_n) \varphi_\infty, \pi_\infty (b_n) \varphi_\infty\rangle\\
&\longrightarrow \langle \pi_0 (Y^p dU (f_l) Y^p) \pi_0 (Y)^{-p} \psi, \pi_0 (Y)^{-p} \psi\rangle = \langle dV (f_l) \psi, \psi\rangle,
\end{align*}
where the last equality follows from equation (\ref{pib}). Since $b_n^\ast Y^p dU(f_l) Y^p b_n \in \cC_f$ and hence $F(b^\ast_n Y^p dU(f_l) Y^p b_n)\ge 0$, we get $\langle dV(f_l)\psi, \psi\rangle  \ge 0$.

Since $\pi_F$ and hence $\pi_0$ are cyclic representations, $\Hh_0$ is separable. Since the Lie group $G$ is connected, $G$ is separable. Therefore, the unitary representation $V$ on $\Hh_0$ can be decompered as a direct integral $\int^\oplus_\Lambda U_\lambda d\mu (\lambda)$ of irreducible unitary representations (see \cite{Ki}, p.127). Now we need two (known) technical results from decomposition theory (see e.g. \cite{S2}, Chapter 12, pp. 343-344 and \cite{Nu}). The first one states that
\begin{equation}\label{dvint}
dV = \int^\oplus_\Lambda dU_\lambda d \mu (\lambda).
\end{equation}
For a unitary representation $W$ of $G$, let $\tau_W$ denote the metric locally convex topology on $\cD^\infty (W)$ given by family of seminorms $\| dW(x^{n_1}_1 {\cdots}x^{n_d}_d) \cdot\|$, where $n_1,{\dots},$ n$_d \in \dN_0$. Since $\Hh_0$ and hence $\cD^\infty (V)$ is separable, there is a countable dense subset $\{\eta_n; n \in \dN\}$ of $\cD^\infty (V) [\tau_V]$. The second result states that  then $\{\eta_n (\lambda); n \in \dN\}$ is dense in $\cD^\infty (U_\lambda) [\tau_{U_\lambda}]~ \mu$- a.e. 

\mn
For $\zeta \in L^\infty (\Lambda, \mu)$, let $M_\zeta$ denote the associated diagonalisable operator on \break
$\Hh_0 = \int^\oplus_\Lambda \Hh (U_\lambda) d\mu (\lambda)$. From (\ref{dvint}) we obtain
$$
\langle dV (f_l) M_\zeta \eta_n, M_\zeta  \eta_n \rangle = \int_\Lambda |\zeta (\lambda)|^2 \langle dU_\lambda (f_l) \eta_n (\lambda), \eta_n (\lambda)\rangle d \mu (\lambda)
$$
for all $\zeta \in L^\infty (\Lambda, \mu)$. Since $ dV (f_l) \ge 0$, the latter implies that there is a $\mu$-null set $N_0$ of $\Lambda$ such that $\langle dU_\lambda (f_l) \eta_n (\lambda), \eta_n (\lambda)\rangle\ge 0$ for all $n\in \dN$ and  all $\lambda \in \Lambda \backslash N_0$. From the density of  $\{\eta_n (\lambda); n \in \dN\}$ in $\cD^\infty (U_\lambda)[\tau_{U_\lambda}]\  \mu$-a.e. it follows that there is a $\mu$-null set $N$ such that $d U_\lambda (f_l) \ge 0$ for all $\lambda \in \Lambda\backslash N$ and $l = 1, {\dots}, r$. That is, the equivalence class of $U_\lambda$ is in  $\cK_f$ for  $\lambda \in \Lambda\backslash N$. Since $d U_\alpha (c-\varepsilon\cdot 1) \ge 0$ for $\alpha \in \cK_f$ by assumption (i), from the latter and (\ref{dvint}) we conclude that $dV (c -\varepsilon \cdot 1) \ge 0$. Therefore, by (\ref{pib}) we have
$$
\langle \pi_0 (Y^n dU (c) Y^n) \psi,\psi \rangle = \langle dV (c) \pi_0 (Y)^n \psi, \pi_0 (Y)^n \psi \rangle
\ge \varepsilon \| \pi_0 (Y)^n \psi\|^2
$$
for $\psi \in \cD^\infty (V)$ and hence
\begin{equation}\label{pic1}
\langle \pi_0 (Y^n dU (c) Y^n) \varphi_0, \varphi_0\rangle \ge \varepsilon \| \pi_0 (Y)^n \varphi_0 \|^2.
\end{equation}

Next we consider the second summand in (\ref{fpix}) for  $x = Y^n dU(c) Y^n$. Let $E(\cdot )$ denote the spectral measure $E^+ (\cdot) \oplus E^- (\cdot)$ on $\Hh_\infty = \Hh^+_\infty \oplus \Hh^-_\infty$. Since $\pi_\infty (\cX_0) = \{0\}$, it follows from formulas (\ref{r8}), (\ref{pinullinf}) and (\ref{piyinf}) that
\begin{align}\label{pic2}
&\langle \pi_\infty (Y^n dU (c) Y^n) \varphi_\infty, \varphi_\infty \rangle =\langle \pi_\infty (Y^n dU (c_{4n}) Y^n)  \varphi_\infty, \varphi_\infty \rangle \nonumber \\
& =\int_{S^d} c_{4n} (t) d \langle E (t) \varphi_\infty, \varphi_\infty\rangle.
\end{align}
By assumption (ii), $c_{4n} (t) > 0$ for all $t \in S^d$. Since $F (Y^n dU(c) Y^n)\le 0$, we conclude from (\ref{fpix}), (\ref{pic1}) and (\ref{pic2}) that $\pi_0 (Y)^n \varphi_0 = 0$ and $\langle E (\cdot) \varphi_\infty, \varphi_\infty \rangle
= 0$. Therefore, $\varphi_0 = 0$ and $\varphi_\infty = 0$, so that $F\equiv 0$ by (\ref{fpix}). Since $F$ is  a state on $\cX$, we have a contradiction. \hfill $\Box$
\begin{thl}\label{L7}
Let $n \in \dN_0$. For arbitrary elements $z_1,{\dots},z_q \in \cX$ there exists $\gs \in \cS := dU(S)$ such that $z_1 A^n_0\gs, {\dots}, z_q A^n\gs \in \cU := dU (\cE (g))$.
\end{thl}
{\bf Proof.} First we prove the assertion for single elements $z \in \cX$ of the form $Y_{j_1l_1}{\cdots} Y_{j_kl_k}$. 
We use induction on $k$. Suppose that the assertion holds for $k$. Let $z = Y_{jl} w$, where $w = Y_{j_1l_1} {\cdots} Y_{j_kl_k}$. By the induction hypothesis, there is $\gs^\prime \in \cS$ such that $w A^n\gs^\prime \in \cU$. Assume that $j\le 0, l\le 0$. Since $S$ is a right Ore set containing $a$ and $dU(a)=A$, there are elements $\gs^{\prime\prime} \in \cS$ and $v \in \cU$ such that $Av = X_jX_l (w A^n\gs^\prime) \gs^{\prime\prime}  $. Set $\gs = \gs^\prime \gs^{\prime\prime} $. Then $z A^n \gs =A^{-1} X_j X_l w A^n \gs^\prime \gs^{\prime\prime} = A^{-1} Av = v \in \cU$. The case $j>0,l>0$ and the case $k=1$ are treated similarly.

It suffices to prove the assertion of Lemma \ref{L7} for element $z_j$ of the form $Y_{j_1l_1}{\cdots} Y_{j_kl_k}$ because these elements and $I$ span $\cX$. We proceed by induction on $q$. For $q=1$ the assertion is proved in the preceding paragraph. Suppose that the assertion is true for $q$. Let $z_1, {\dots}, z_{q+1} \in \cX$. Then there exist elements $\gs_1,\gs_2 \in \cS$ such that $z_l A^n \gs_1 \in \cU$ for $l=1,{\dots},q$ and $z_{q+1} A^n \gs_2 \in \cU$. By the right Ore property of $S$, there are $\gs_3 \in \cS$ and $u \in \cU$ such that $\gs_1\gs_3 = \gs_2 u=:\gs$. Then, $\gs \in \cS, z_l A^n \gs = (z_l A^n \gs_1)\gs_3 \in \cU$ for $l=1,{\dots},q$ and $z_{q+1} A^n\gs =(z_{q+1} A^n \gs_2) u \in \cU$.\hfill $\Box$

\mn
Now we are able to complete the proof of Theorem 1.1. By Lemma \ref{L6}, there exist finitely many elements $z_{jl} \in \cX, l=0,{\dots},d$, such that
$$
Y^n dU(c)Y^n = \sum_j z^\ast_{j0} z_{j0} + \sum_j \sum^d_{l=1} z^\ast_{jl} Y^p dU(f_l)Y^p z_{jl}.
$$
Let $\gs = dU(s), s\in \cS$. Multiplying both sides by $A^n\gs$ from the right and by $(A^n\gs)^\ast$ from the left we obtain 
$$
dU (s^\ast c s) = \sum_j (z_{j_0} A^n \gs)^\ast z_{j_0} A^n \gs + \sum_j \sum^d_{l=1} (Y^p z_{jl} A^n \gs)^\ast dU(f_l) Y^p z_{jl} A^n\gs.
$$

By Lemma \ref{L7} we can find $\gs = dU(s)$ such that $z_{j0} A^n\gs \in dU(\cE(g))$ and $Y^p P z_{jl} A^n\gs \in dU(\cE(g))$ for all $j,l$. Then the right-hand side of the preceding equation is in $dU(\cT_f)$. Now we choose the unitary representation $U$ of $G$ such that the representation $dU$ of $\cE(g)$ is faithful (for instance, it suffices to take the regular representation of $G$). Then it follows from $dU(s^\ast c s)\in dU(\cT_f)$ that $s^\ast c s \in \cT_f$.

Finally, we suppose that $m$ is odd. Then $c^\prime := \sum^d_{k=0} x^\ast_k cx_k$ satisfy assumptions (i) and (ii) with $m^\prime = m+1$ even, so the assertion of Theorem 1.1 follows from the previous case. \hfill $\Box$
\mn


\begin{thebibliography}{9999} 


\bibitem[D]{D}
Dixmier, J., 
{\it Alg\'ebres Enveloppantes}, Gauthier-Villars, Paris, 1974.

\bibitem[Ki]{Ki}
Kirillov, A.A., 
{\it Elements of the Theory of Representations}, Springer-Verlag, Berlin, 1976.


\bibitem[K]{K}
K\"othe, G.,
{\it Topological Vector Spaces II}, Springer-Verlag, Berlin, 1979.


\bibitem[M]{M}
Marshall, M.,
{\it Positive Polynomials and Sums of Squares}, Univ. Pisa, Dipart. Mat., Istituti Editoriali e Poligrafici Internaz., 2000.

\bibitem[N]{N}
Nelson, E.,
{\it Analytic vectors}, Ann. Math. {\bf 70} (1959), 572 -- 615.

\bibitem[Nu]{Nu}
Nussbaum, A. E.,
{\it Reduction theory for unbounded closed operators in Hilbert space}, Duke Math. {\bf 31} (1964), 33 -- 44.

\bibitem[NS]{NS}
Nelson, E. and W. F. Stinespring,
{\it Representation of elliptic operators in an enveloping algebra}, Amer. J. Math. {\bf 81} (1959), 547 -- 560.

\bibitem[P]{P}
Putinar, M.,
{\it Positive polynomials on compact semi-algebra sets}, Indiana Univ. Math. J. {\bf 42} (1993), 969 -- 984.

\bibitem[PD]{PD}
Prestel, A. and C.N. Delzell,
{\it Positive Polynomials}, Springer-Verlag, Berlin, 2001.

\bibitem[PV]{PV}
Putinar, M. and F.-H. Vasilescu,
{\it Solving moment problems by dimensional extension}, Ann. Math. {\bf 149} (1999), 1087 -- 1107.

\bibitem[R]{R}
Reznick, B.,
{\it Uniform denominators  in Hibert's Seventeenth problem}, Math. Z. {\bf 220} (1995), 75 -- 97.

\bibitem[S]{S}
Scheiderer, K.,
{\it Positivity and sums of squares: A guide to some recent results.} Preprint, Duisburg, 2002.

\bibitem[S1]{S1}
Schm\"udgen, K.,
{\it The $K$-moment problem for compact semi-algebraic sets}, Math. Ann. {\bf 289} (1991), 203 -- 206.

\bibitem[S2]{S2}
Schm\"udgen, K.,
{\it Unbounded Operator Algebras and Representation\\ Theory}, Birkh\"auser-Verlag, Basel, 1990.

\bibitem[S3]{S3}
Schm\"udgen, K.,
{\it A strict Positivstellensatz for the Weyl algebra}, Preprint, Leipzig, 2004, math.AC/0403076.

\bibitem[Sz]{Sz}
Szymanski, W.,
{\it Group $C^\ast$-algebras as algebras of "continuous functions" with non-commuting variables}, Proc. Amer. Math. Soc. {\bf 107} (1989), 353 -- 359.

\bibitem[Wa]{Wa}
Warner, G.,
{\it Harmonic analysis on semi-simple Lie groups I}, Springer-Verlag, Berlin, 1972.

\end{thebibliography}
\end{document}